\documentclass[12pt,reqno,intlim]{amsart}

\usepackage{amssymb,amsmath}

\usepackage[arrow,matrix,curve]{xy}

\newdir{ >}{{}*!/-5pt/\dir{>}}

\newcommand{\nc}{\newcommand}

\nc{\bC}{\bold{C}} \nc{\bN}{\Bbb{N}} \nc{\cF}{\mathcal{F}}
\nc{\cE}{\mathcal{E}} \nc{\cR}{\mathcal{R}} \nc{\cM}{\mathcal{M}}
\nc{\al}{\alpha} \nc{\bt}{\beta} \nc{\gm}{\gamma} \nc{\dl}{\delta}
\nc{\om}{\omega} \nc{\sg}{\sigma} \nc{\Sg}{\Sigma} \nc{\vf}{\varphi}
\nc{\ve}{\varepsilon} \nc{\os}{\overset} \nc{\ol}{\overline}
\nc{\ul}{\underline} \nc{\us}{\underset} \nc{\sbs}{\subset}
\nc{\bsl}{\backslash} \nc{\Ra}{\Rightarrow}
\nc{\lra}{\longrightarrow} \nc{\all}{\allowdisplaybreaks}

\nc{\Codes}{\operatorname{{\bold{Codes}}}}
\nc{\RegMono}{\operatorname{\mathcal{R}{\rm{eg}\mathcal{M}{\rm{ono}\!}}}}
\nc{\RegEpi}{\operatorname{\mathcal{R}{\rm{eg}\mathcal{E}{\rm{pi}\!}}}}
\nc{\Mn}{\operatorname{\mathcal{M}{\rm{ono}\!}}}
\nc{\Ep}{\operatorname{\mathcal{E}{\rm{pi}\!}}}
\nc{\Rg}{\operatorname{\mathcal{R}{\rm{eg}\!}}}
\nc{\Ob}{\operatorname{Ob\!}}

\numberwithin{equation}{section}

\newtheorem{theo}{\ \ \ Theorem}[section]
\newtheorem{lem}[theo]{\ \ \ Lemma}
\newtheorem{prop}[theo]{\ \ \ Proposition}

\theoremstyle{definition}

\theoremstyle{remark}

\begin{document}

\title[]
{Descent in the dual category of ternary rings}

\author{Guram Samsonadze and Dali Zangurashvili}

\maketitle

% Abstract.
\begin{abstract}
It is shown that, in the variety of ternary rings, the elements of amalgamated free products have unique normal forms, and, moreover, this variety satisfies the strong amalgamation property. Applying these statements, effective codescent morphisms of ternary rings are characterized. In view of the fact that the category of ternary rings contains the category of commutative associative unitary rings as a full subcategory, the class of effective codescent morphisms in the latter category (which, according to the well-known Joyal-Tierney's criterion, are precisely monomorphisms $R\rightarrowtail S$ which are pure as monomorphisms of $R$-modules) is compared with that of morphisms between commutative associative unitary rings which are effective codescent in the category of ternary rings. It turns out that the former class is contained in the latter one, but does not coincide with it. 
 
\bigskip


\end{abstract}
\vskip+2mm

% 1.
\section{Introduction}

The aim of the present paper is to characterize effective descent morphisms in the dual category of ternary rings. We use the term `ternary ring' in the sense of \cite{BH} by Chajda and Hala\v{s}, which slightly differs from that of \cite{M1}, \cite{M2}, \cite{M3} by Machala, where the notion of a ternary ring was introduced and applied for the coordinatization of projective and affine planes with homomorphism  (but the roots of this notion go back to the earlier work \cite{H} by Hall). Chajda and Hala\v{s} define a ternary ring\footnote{The unique difference of this definition with that by Machala is that Machala requires, in addition to the conditions (i) and (ii) below, that  $0\neq 1$.} as a set $R$ equipped with a ternary operation $t$ and two constants $0$ and $1$ such that the following conditions hold:\vskip+1mm

(i) the identities 
\begin{equation}
t(0,x,y)=y,
\end{equation}
\begin{equation}
t(x,0,y)=y,
\end{equation}
\begin{equation}
t(1,x,0)=x,
\end{equation}
\noindent and
\begin{equation}
t(x,1,0)=x,
\end{equation}
\noindent are satisfied;
\vskip+1mm
(ii) for any $a,b,c\in R$ there exists a unique $x\in R$ with 
\begin{equation}
t(a,b,x)=c.
\end{equation}
\vskip+2mm
We see that a ternary ring is not a ring in the traditional sense. However the category of ternary rings contains the category of traditional rings with unit as a full reflective subcategory: if $R$ is a ring of the latter kind, then the same set $R$ equipped with the operation  $$t(a,b,c)=ab+c$$
\noindent is a ternary ring \cite{M1}.
%(see, e.g., the paper \cite{BH} by Chajda and Hala\v{s}).

 As noted by Chajda and Hala\v{s} \cite{BH}, the category of ternary rings  is a variety of universal algebras\footnote{When referring to a ternary ring as to a universal algebra, the authors of \cite{BH} use the term `bi-ternary ring'.}. The signature $\mathfrak{F}$ of this variety consists of two ternary operations $t$ and $q$, and two nullary operations $0$ and $1$, while the set $\Sigma$ of identities consists of (1.1)-(1.4) and 
\begin{equation}
q(x,y,t(x,y,z))=z,
\end{equation}
\begin{equation}
t(x,y,q(x,y,z))=z.
\end{equation}

\vskip+2mm

The present paper continues the series of works  \cite{Z1}, \cite{Z3}, \cite{Z4}, \cite{SZ} on the descent problem (i.e., the problem of characterizing effective descent morphisms) in categories dual to varieties of universal algebras posed by Janelidze and motivated by the classical results on effective descent morphisms in the category dual to the variety of commutative associative unitary rings (referred to as `commutative rings', in the sequel). Recall that the first result on this issue, still used in commutative algebra and algebraic geometry, was obtained by Grothendieck. This result asserts that a monomorphism $p:R\rightarrowtail S$ of such rings is an effective codescent morphism (i.e., an effective descent morphism in the dual category) if $p$ makes $S$ a faithfully flat $R$-module. The complete solution to the descent problem for the dual of this variety was given by Joyal and Tierney. They found that a monomorphism $p$ of commutative rings is an effective codescent morphism if and only if it is a pure monomorphism of $R$-modules  (unpublished). Later Mesablishvili gave a new proof of this result 
\cite{M} (see also the paper \cite{JT2} by Janelidze and Tholen).

In \cite{Z1}, \cite{Z3}, \cite{Z4}, Zangurashvili dealt with abstract and particular varieties satisfying the amalgamation property (note that this property is not satisfied by the variety of commutative rings, but is satisfied by the variety of ternary rings, as is shown in this paper). In particular, she characterized effective codescent morphisms of groups, loops, (left/right-)quasigroups, and magmas. Employing results of the papers \cite{Z1}-\cite{Z3}, in our joint work \cite{SZ}, we related the descent problem for categories dual to varieties satisfying the amalgamation property to the notion of a confluency, which is one of the central notions in term rewriting theory.

 In the present paper, applying the techniques developed in \cite{SZ}, we prove that every codescent morphism of ternary rings is effective. To this end, we show for that variety that the elements of amalgamated free products have unique normal forms. This fact, in particular, implies that the variety of ternary rings satisfies the strong amalgamation property. These results combined with those by Zangurashvili \cite{Z1} and by Chajda and Halaš \cite{BH} on (universal-algebraic) 0-ideals of ternary rings leads to the characterization of effective codescent morphisms in the category of ternary rings. It asserts that a monomorphism $R\rightarrowtail S$ is such a morphism if and only if it satisfies the $0$-ideal extension property.

Finally, in view of the fact that the category of ternary rings contains the category of commutative rings as a full subcategory, we compare the class of effective codescent morphisms in the latter category  with that of morphisms between commutative rings which are effective codescent in the category of ternary rings. It turns out that the latter class coincides
 with the class of commutative ring monomorphisms that satisfy the traditional ring-theoretic ideal extension property. Note that this class (i.e., the class of monomorphisms which satisfy this property) naturally appeared in the study of the class of commutative ring monomorphisms which satisfy the above-mentioned Joyal-Tierney's criterion, from a purely ring-theoretic point of view decades ago. Since then, the question whether these two classes are equal has been studied by a number of authors. Applying their results, we find that if a monomorphism is an effective codescent morphism in the category of commutative  rings, then it is such also in the category of ternary rings. However, the converse is not true.

The second author gratefully acknowledges the financial support of Shota Rustaveli National Science Foundation of Georgia (FR-22-4923).

\section{Preliminaries}
We assume that the reader is familiar with the basics of term rewriting theory (see, e.g., the book \cite{BN} by Baader and Nipkow). 

 Before continue, let us agree on the notation. Let $\mathfrak{F}$ be a signature. Let $X$ be a countable set with $\mathfrak{F} \cap X=\varnothing$. For  terms $\tau$, $\tau'$ over $X$ and a position $p$ in $\tau$, the symbol $\tau\mid_p$ denotes the subterm on the position $p$ of $\tau$, while the symbol $\tau[\tau']_p$ denotes the term obtained by replacing the subterm  $\tau\mid_p$ by the term $\tau'$ in $\tau$. The symbol $Var(\tau)$ denotes the set of variables which present in the term $\tau$.

 Assume that $\Sigma$ is a set of ($\mathfrak{F}$-)identities with variables from $X$, and $\mathcal{C}$ is the variety of universal algebras determined by $\mathfrak{F}$ and $\Sigma$. 

From now on all identities are assumed to be oriented.

Let $\mathfrak{F}_0$ denote the set of all nullary operations from $\mathfrak{F}$. In \cite{SZ}, we introduced the following conditions:\vskip+2mm
 
 (*) for any identity $l=r$ from $\Sigma$, no variable occurs in the term $r$  more often than in the term $l$, and moreover, the size of $l$ is greater than the size of $r$;\vskip+2mm
  
 (**) if the set $\mathfrak{F}_0$ is not empty, then, for any nontrivial algebra $A$ from the variety determined by $\mathfrak{F}$ and $\Sigma$, the mapping $\mathfrak{F}_0\rightarrow A$ sending a nullary operation to its value in $A$ is injective;\vskip+2mm
 
 (***) for any identity $l=r$ from $\Sigma$, any subterm $l'$ of $l$ that is neither a variable nor a nullary operation symbol, we have $Var(l')=Var(l)$.
\vskip+4mm

Let $I$ be a nonempty set, and $A_i$ be a $\mathcal{C}$-algebra, for any $i\in I$. Let $B$ be a subalgebra of all $A_i$ with the embeddings $m_i:B\rightarrowtail A_i$. Without loss of generality we can assume that $A_i\cap X=A_i\cap \mathfrak{F}=\emptyset$, and $A_i\cap A_j=B$, for distinct $i$ and $j$.

On the set of terms over the set $\underset{i\in I}\cup A_i$ we introduce the binary relation $\rightsquigarrow$ as follows. We say that $\tau\rightsquigarrow \tau'$ if one of the following conditions  is satisfied:\vskip+2mm
($\textbf{C}$) there exists an identity $l=r$ in $\Sigma$, a substitution $\sigma$ (on the set $\underset{i\in I}\cup A_i$) and a position $p$ of $\tau$ such that $\tau\mid_p=\sigma(l)$ and $\tau'=\tau[\sigma(r)]_p$;
\vskip+2mm

($\textbf{C1}$) we have $\tau\mid_p=o(a_1,a_2,...,a_n)$, for some position $p$ in $\tau$, an $n$-ary operation symbol $o$ from the signature $\mathfrak{F}$, and some elements $a_1,a_2,$ $...,a_n$ from one and the same algebra $A_i$, and, moreover, $\tau'=\tau[a]_p$, where $a$ is the value of the term $o(a_1,a_2,...,a_n)$ in $A_i$. \vskip+3mm

 Let $A$ be the free product of algebras $A_i$ with the amalgamated subalgebra $B$, and let $n_i: A_i\rightarrow A$ be the canonical homomorphisms, for any $i\in I$. It is well-known that $A$ is isomorphic to the quotient of (not necessarily $\mathcal{C}$-)algebra of terms over the set $\underset{i\in I}\cup A_i$ with respect to the reflexive transitive closure of the relation $\rightsquigarrow$. The condition (*) implies that any element $\alpha$ of $A$ can be written as a term $\tau$ over the set $\underset{i\in I}\cup A_i$ such that there is no a term $\tau'$ over the same set with $\tau\rightsquigarrow \tau'$. The term $\tau$ is called a \textit{normal form} of the element $\alpha$. The question that arises is whether normal forms are unique.\vskip+2mm

\begin{theo}\cite{SZ}
 Let $\mathcal{C}$ be a variety of universal algebras represented by a confluent term rewriting system $(\mathfrak{F},\Sigma)$ satisfying the conditions (*)-(***). Then elements of amalgamated free products have unique normal forms in $\mathcal{C}$.
 \end{theo}

A variety $\mathcal{C}$ is said to satisfy the strong amalgamation property if, for any family $(A_i)_{i\in I}$ of $\mathcal{C}$-algebras, the homomorphisms $n_i$'s are monomorphisms, and $n_i(A_i)\cap n_j(A_j)=n_im_i(B)$, for any distinct $i$ and $j$. This definition is equivalent to the one where only the families  $(A_i)_{i\in I}$ with $card(I)=2$ are considered (see, e.g., the paper \cite{I} by Imaoka).
\vskip+2mm 
 
  \begin{theo} \cite{SZ}
  Under the conditions of Theorem 2.1, the variety $\mathcal{C}$ satisfies the strong amalgamation property.
 \end{theo}

\vskip+5mm
 Recall the definition of an (effective) codescent morphism \cite{JT}. Let $\mathcal{C}$ be a category with pushouts, and let $p:B\rightarrow E$ be its morphism. There is the so-called change-of-cobase functor 
\begin{equation}
 p_{*}:B/\mathcal{C}\rightarrow E/\mathcal{C},
\end{equation}
\noindent which sends a $\mathcal{C}$-morphism $\varphi:B\rightarrow C$ to the pushout of $\varphi$ along $p$. It is well-known that the functor $p_{*}$ has a right adjoint $p^{!}$. It sends a morphism $\psi:E\rightarrow C$ to the composition $\psi p$. A morphism $p$ is called a \textit{codescent morphism} if the functor $p_{*}$ is premonadic, and is called an \textit{effective codescent morphism} if $p_{*}$ is monadic.

The following theorem by Janelidze and Tholen is well-known.
 
\begin{theo}\cite{JT1}, \cite{JT}
Let $\mathcal{C}$ be a category with pushouts and equalizers. A morphism $p$ is a codescent morphism if and only if it is a couniversal regular monomorphism. i.e. a morphism whose any pushout is a regular monomorphism.
\end{theo}

For the definitions of an ideal of an algebra from a variety of universal algebra and of an ideal-determined variety, we refer the reader to the paper \cite{GU} by Gumm and Ursini.
\vskip+1mm
 
Applying Theorem 2.3, Zangurashvili proved the following

 \begin{theo} \cite{Z1}
 Let $\mathcal{C}$ be a variety of universal algebras\footnote{In fact, the paper \cite{Z1} deals with the general case of an arbitrary category with pushouts and equalizers.} that satisfies the strong amalgamation property. For a monomorphism $p:B\rightarrowtail E$, the conditions (i) and (ii) below are equivalent.  If $\mathcal{C}$ has a nullary operation $0$ and is $0$-ideal determined, then they are equivalent also to the condition (iii):\vskip+2mm
 
 (i) the monomorphism $p$ is a codescent morphism;\vskip+2mm
 
 (ii) the monomorphism $p$ satisfies the congruence extension property, i.e., for any congruence $\theta$ on $B$, there is a congruence $\theta'$ on $E$ such that $\theta'\cap (B\times B)=\theta$.\vskip+2mm

 (iii) the monomorphism $p$ satisfies the $0$-ideal extension property, i.e., for any $0$-ideal $I$ on $B$, there is a $0$-ideal $I'$ on $E$ such that $I'\cap B=I$.\vskip+2mm
 \end{theo}

In \cite{Z3}, Zangurashvili proved that if, in a variety of universal algebras, elements of amalgamated free products have unique normal forms, then all codescent morphisms are effective in it. This enabled us to formulate the following
 \begin{theo}\cite{SZ}
 Let $\mathcal{C}$ be a variety of universal algebras represented by a confluent term rewriting system $(\mathfrak{F},\Sigma)$ satisfying the conditions (*)-(***). Then every codescent morphism of $\mathcal{C}$ is effective.
 \end{theo}

\section{The variety of ternary rings}
Let $\mathcal{C}$ now be the variety of ternary rings. Chajda and Hala\v{s} noted that the following identities are satisfied in this variety \cite{BH}:
\begin{equation}
q(0,x,y)=y,
\end{equation}
\begin{equation}
q(x,0,y)=y,
\end{equation}
\begin{equation}
q(1,x,x)=0,
\end{equation}
\begin{equation}
q(x,1,x)=0.
\end{equation}
\vskip+2mm
\noindent One can verify that these are precisely the identities corresponding to all non-joinable critical pairs with respect to the signature $\mathfrak{F}$ defined in Section 1 and the set $\Sigma$ of oriented identities (1.1)-(1.4), (1.6), (1.7). Adding (3.1)-(3.4) to $\Sigma$, we still obtain a set $\Sigma'$ of identities that determines the variety of ternary rings. The condition (*) from Section 2 is obviously satisfied by $\Sigma'$, and hence the term rewriting system $(\mathfrak{F}, \Sigma')$ is terminating. One can verify that all its critical pairs are joinable. Thus, applying the Critical Pair Theorem (see, e.g.,\cite{BN}), we obtain

\begin{lem}
The term rewriting system $(\mathfrak{F},\Sigma')$ is confluent.
\end{lem}

The term rewriting system $(\mathfrak{F},\Sigma')$ obviously satisfies also the conditions (**) and (***). Therefore, applying Theorems 2.1, 2.2 and Lemma 3.1 we obtain 
 
 \begin{theo}
  In the variety of ternary rings, the elements of amalgamated free products have unique normal forms (with respect to the representation $(\mathfrak{F},\Sigma')$).
  \end{theo}
 
 \begin{theo}
 The variety of ternary rings satisfies the strong amalgamation property.
 \end{theo}
 Recall the following result by Chajda and Hala\v{s}.
 \begin{theo}\cite{BH}
 The variety of ternary rings is $0$-ideal determined.
 \end{theo}

Chajda and Hala\v{s} characterized $0$-ideals of ternary rings. To this end, they introduced the following terms:
$$\tau_1(x,y)=q(1,t(1,x,y),x),$$
$$\tau_2(x_1,x_2,x_3,y)=q(1,x_3,q(x_1,x_2,t(1,t(x_1,x_2,x_3),y))),$$
$$\tau_3(x_1,x_2,x_3,y_1,y_2,y_3)=$$
$$q(1,t(x_1,x_2,x_3),t(t(1,x_1,y_1),t(1,x_2,y_2),t(1,x_3,y_3))).$$
%\noindent where the symbol $x\circ y$ denotes the term $t(1,x,y)$. 

\begin{theo} \cite{BH}
 A nonempty subset $I$ of a ternary ring $R$ is a $0$-ideal if and only if, for any $a_1,a_2,a_3 \in R$ and any $h_1,h_2,h_3\in I$, we have $\tau_1(a_1,h_1),\tau_2(a_1,a_2,a_3,h_1),\tau_3(a_1,a_2,a_3,h_1,h_2,h_3)\in I$. 
 \end{theo} 
% \begin{equation}
 \vskip+2mm
Theorem 3.5 implies that, for any ternary ring $S$ and its any non-empty subset  $Y$, the smallest ideal containing $Y$ obviously is the union $\underset{i\geq 0} \cup Y_{i}$, where $Y_0=Y$ and 
$$Y_i=Y_{i-1}\cup \lbrace{\tau_1(a_1,h_1),\tau_2(a_1,a_2,a_3,h_1),\tau_3(a_1,a_2,a_3,h_1,h_2,h_3):}$$
\begin{equation}
{a_1,a_2,a_3\in S, h_1,h_2,h_3\in Y_{i-1}}\rbrace,
\end{equation}

\noindent for any  $i\geq 1$. Hence, Theorems 2.4, 2.5 and Lemma 3.1 imply
 \begin{theo}
Let $R$ and $S$ be ternary rings, and $p:R\rightarrowtail S$ be a monomorphism in $\mathcal{C}$. The following conditions are equivalent:\vskip+2mm

(i) $p$ is an effective codescent morphism;\vskip+2mm

(ii) $p$ is a codescent morphism;\vskip+2mm

(iii) $p$ satisfies the congruence extension property;\vskip+2mm

(iv) $p$ satisfies the $0$-ideal extension property;

\vskip+2mm

(v) for any $0$-ideal $I$ of $R$, we have $(\underset{i\geq 0} \cup Y_{i}) \cap R=I$, where $Y_0=I$ and, for all $i\geq 1$, the sets $Y_i$ are defined by equality (3.5).
\end{theo}

\section{The case of commutative rings}

There naturally arises the question: what relationship exists between the class of effective codescent morphisms in the variety of commutative rings and that of monomorphisms of commutative rings which are effective codescent morphisms in the variety of ternary rings? To answer this question, we will not apply the known categorical results on the reflection of effective codescent morphisms by a functor since all of these results require that the functor preserves all or at least some pushouts \cite{JT}, \cite{R}, \cite{M4}, \cite{Z2}. At that, the embedding functor from the category of commutative rings to the category of ternary rings does not preserve pushouts even of monomorphisms (otherwise the category of commutative rings would satisfy the amalgamation property, but it does not). For that reason, we will deal with this problem via compering the class of monomorphisms determined by the condition (iv) of Theorem 3.6 with the class of monomorphisms from the Joyal-Tierney's criterion (see Introduction).

For a ring $R$ (in the traditional sense) with unit, the terms $\tau_1$, $\tau_2$ and $\tau_3$ take the forms:
$$\tau_1(x,y)=-y;$$
$$\tau_2(x_1,x_2,x_3,y)=y;$$
$$\tau_3(x_1,x_2,x_3,y_1,y_2,y_3)=x_1y_2+y_1x_2+y_1y_2+y_3.$$
\vskip+1mm
\noindent This immediately implies 

\begin{lem} Let $R$ be a (traditional) associative ring. A subset $I$ of $R$ with $0\in I$ is a $0$-ideal of $R$ as of a ternary ring if and only if $I$ is a two-sided ideal of $R$ as of a traditional ring.
\end{lem}

 Therefore, the term `$0$-ideal' in the condition (iv) of Theorem 3.6 can be understood in the ring-theoretic context. Hence we have
 \begin{lem}
Let $R$ and $S$ be commutative rings. A monomorphism $R\rightarrowtail S$ is an effective codescent morphism in the variety of ternary rings if and only if it satisfies the traditional ring-theoretic ideal extension property. 
 \end{lem} 
 
The characterization of pure monomorphisms of modules over commutative rings in terms of solvability of some systems of linear equations by Cohn is well-known (see, e.g., \cite{L}). It implies that a monomorphism $R\rightarrowtail S$ of commutative rings that is pure as a homomorphism of $R$-modules satisfies the (traditional ring-theoretic) ideal extension property. However the converse is not true, as it follows from the counterexamples constructed by Enochs (given in \cite{GM}), Lady (given in \cite{B}), and Emmanouil (given in \cite{L}). But the converse is true if $R$ is a principal ideal domain \cite{L}. In this statement, ``$R$ is a principal ideal domain" can be replaced by the condition that $S$ is projective as an $R$-module, according to the result by Fieldhouse given in \cite{L}. Considering these facts and Theorem 3.6, we obtain the following

\begin{prop}
Let $R$ and $S$ be commutative rings, and $R\rightarrow S$ be a homomorphism. If it is an effective codescent morphism in the category of commutative rings, then it is such also in the category of ternary rings. The converse is not true, in general. But it is true if either $R$ is a principal ideal domain or $S$ is projective as an $R$-module.
\end{prop}

Finally note that the above-mentioned results imply that the condition that $\mathcal{C}$ satisfies the strong amalgamation property in Theorem 2.4 (and hence in Proposition 2.2 of \cite{Z1}) is essential.

\vskip
+3mm

 \vskip+2mm
 
\textit{Authors' addresses:}
 
\textit{Guram Samsonadze, Georgian Technical University, 77 Kostava Str., Tbilisi, 0160, Georgia; e-mail: g.samsonadze@gtu.ge}
\vskip+1mm
\textit{Dali Zangurashvili, A. Razmadze Mathematical Institute of Tbilisi State University, 6 Alexidze Str., 0193, Georgia;}

\textit{e-mail: dali.zangurashvili@tsu.ge}

\end{document}